\documentclass[12pt,reqno]{amsart}

\usepackage[usenames]{color}
\usepackage{amssymb}
\usepackage{amsmath}
\usepackage{amsthm}
\usepackage{amsfonts}
\usepackage{amscd}
\usepackage{graphicx}

\usepackage{tikz}
 
\usepackage{color}
\usepackage{fullpage}
\usepackage{float}

\usepackage{graphics}
\usepackage{latexsym}
 
\theoremstyle{plain}
\newtheorem{theorem}{Theorem}

\theoremstyle{definition}

\theoremstyle{remark}

\begin{document}

\title{Partial Dyck paths with Air Pockets}

\author[H.~Prodinger]{Helmut Prodinger}

\address{Helmut Prodinger,
	Mathematics Department, Stellenbosch University,
	7602 Stellenbosch, South Africa, and NITheCS (National Institute for
	Theoretical and Computational Sciences), South Africa.}
\email{hproding@sun.ac.za}

	\begin{abstract} 
		Dyck paths with air pockets are obtained from ordinary Dyck paths by compressing maximal
		runs of down-steps into giant down-steps of arbitrary size. Using the kernel method, we
		consider partial Dyck paths with air pockets, both, from left to right and from right to left.
		
		In a last section, the concept is combined with the concept of skew Dyck paths.
			\end{abstract}
 
 \maketitle

\section{Introduction}
 
 In a paper that was posted on valentine's day~\cite{valentine}, Baril et al.\ introduced a new family
 of Dyck-like paths, called \emph{Dyck paths with air pockets}. Many of the usual 
 parameters that one could think of are investigated in this paper. The paths have the usual 
 up-steps $(1,1)$ and down-steps $(1,-k)$ for any $k=1,2,\dots$, but no such down-steps may follow
 each other. Otherwise, they cannot go into negative territory, and must end at the $x$-axis, as usual.
 One could just think about ordinary Dyck paths, and each (maximal) run of down-steps is condensed into one
 (giant) downstep. 
 
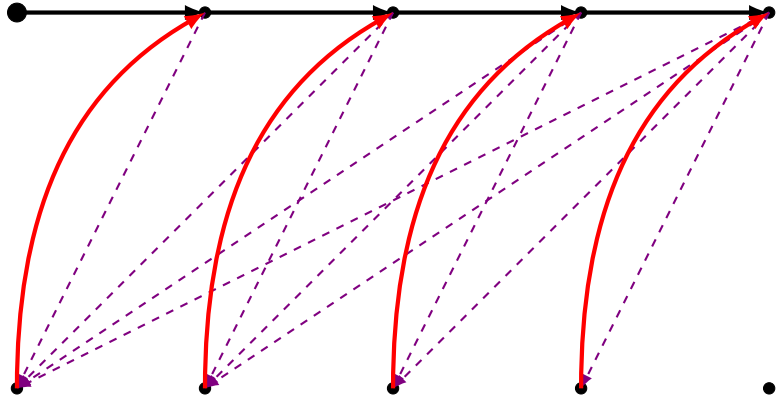
\begin{figure}[ht]\label{myfig}

		\begin{center}
			\begin{tikzpicture}[scale=2.5]
				\draw (0,0) circle (0.05cm);
				\fill (0,0) circle (0.05cm);
				
				\foreach \x in {0,1,2,3,4}
				{
					\draw (\x,0) circle (0.03cm);
					\fill (\x,0) circle (0.03cm);
				}

				\foreach \x in {0,1,2,3,4}
				{
					\draw (\x,-2) circle (0.03cm);
					\fill (\x,-2) circle (0.03cm);
				}
				
				\foreach \x in {0,1,2,3}
				{
					\draw[ ultra thick, -latex] (\x,0) -- (\x+1,0);
					
				}

					\draw[   thick, violet,dashed, -latex] (4,0) -- (3,-2);
					\draw[   thick, violet,dashed, -latex] (4,0) -- (2,-2);
					\draw[   thick, violet,dashed, -latex] (4,0) -- (1,-2);
					\draw[   thick, violet,dashed, -latex] (4,0) -- (0,-2);
					
						\draw[   thick, violet,dashed, -latex] (3,0) -- (2,-2);
						\draw[   thick, violet,dashed, -latex] (3,0) -- (1,-2);
						\draw[   thick, violet,dashed, -latex] (3,0) -- (0,-2);
			 \draw[   thick, violet,dashed, -latex] (2,0) -- (1,-2);
			 \draw[   thick, violet,dashed, -latex] (2,0) -- (0,-2);
			 \draw[   thick, violet,dashed, -latex] (1,0) -- (0,-2);

				  \draw[ ultra   thick, red , -latex] (0,-2) to [out= 90,in= 210] (1,0);
				\draw[ ultra   thick, red , -latex] (1,-2) to [out= 90,in= 210] (2,0);
				 \draw[ ultra   thick, red , -latex] (2,-2) to [out= 90,in= 210] (3,0);
				\draw[ ultra   thick, red , -latex] (3,-2) to [out= 90,in= 210] (4,0);

			\end{tikzpicture}
		\end{center}
		\caption{Graphical description of Dyck paths with air pockets. Top layer describes the situation after an up-step, bottom layer after a down-step. }
	 
	\end{figure}

 The Figure~\ref{myfig} explains the actions readily.
 
	 We introduce generating functions $f_k(z)$ and $g_k(z)$ where  the coefficient of $z^n$ in one of these functions counts paths ending in the respective state according to the number of steps.
	 The function $f_0(z)+g_0(z)$ counts the Dyck paths with air pockets, as the zero in the index just means that they returned to 
	 the $x$-axis.
	 
 In this short paper, we will enumerate \emph{partial} Dyck paths with air pockets, namely we allow the path to end at level $i$.
	 In other words, we compute  all $f_k(z)$ and $g_k(z)$.
	 
	 Our instrument of choice is the kernel method, as can be found in the popular account \cite{prodinger-kernel}.
	 
	 \section{Generating functions}
	 
	 Just looking at Figure \ref{myfig}, we find the following recursion, where we write $f_k$ for $f_k(z)$ for simplicity:
	 \begin{align*}
	 	f_0&=1,\\
	 	f_k&=zf_{k-1}+zg_{k-1},\quad k\ge1,\\
	 	g_k&=zf_{k+1}+zf_{k+2}+zf_{k+3}+\cdots,
	 \end{align*}
 and now we introduce bivariate generating functions 
 \begin{equation*}
 	F(u,z)=F(u)=\sum_{k\ge0}u^kf_k(z), \quad G(u,z)=G(u)=\sum_{k\ge0}u^kg_k(z).
 \end{equation*}
Summing the recursions,
	  \begin{equation*}
F(u)=1+zuF(u)+zuG(u)
	  \end{equation*}
  and
	  \begin{equation*}
	  	G(u)=\sum_{k\ge0}u^kz\sum_{j>k}f_j=z\sum_{j>0}f_j\sum_{k=0}^{j-1}u^k=
	  	z\sum_{j>0}f_j\frac{1-u^k}{1-u}
	   =
  	\frac z{1-u}(F(1)-F(u)) .
  \end{equation*}
Eliminating one function, we are left to analyze
\begin{equation*}
	F(u)=1+zuF(u)+\frac {z^2u}{1-u}(F(1)-F(u)) .
\end{equation*}
Solving, we find
\begin{equation*}
F(u)=\frac {1-u+{z}^{2}uF(1)}{-zu+z{u}^{2}+{z}^{2}u+1-u}
=\frac {1-u+{z}^{2}uF(1)}{z(u-s_1)(u-s_2)},
\end{equation*}
with 
\begin{align*}
	s_1&=\frac {1+z-{z}^{2}+\sqrt {-{z}^{2}-2{z}^{3}-2z+{z}^{4}+1}}
		{2z	},\\
		 s_2&=\frac {1+z-{z}^{2}-\sqrt {-{z}^{2}-2{z}^{3}-2z+{z}^{4}+1}}
		{2z	}.	
\end{align*}
Note that $s_1s_2=\frac1z$. We still need to compute $F(1)$. Before we can plug in $u=1$ and compute it, we must  cancel the bad factor of both, numerator and denominator. In this case, this is the factor $u-s_2$, since the reciprocal of it would not allow a Taylor expansion around $u=1$. The result is
\begin{equation*}
	F(u)= \frac{-1+z^2F(1)}{zs_2-z+z^2-1+zu},	
\end{equation*}
from which we  now can compute $F(1)$ by plugging in $u=1$. We get 
\begin{equation*}
	F(1)= \frac{-1+z^2F(1)}{zs_2+z^2-1}	=\frac{1}{1-zs_2}
\end{equation*}
and therefore
\begin{equation*}
F(u)=\frac{1-s_1}{1-zs_2}\frac{1}{u-s_1}=-\frac{1}{s_1} \frac{1-s_1}{1-zs_2}\frac{1}{1-u/s_1}.
\end{equation*}
Reading off the coefficient of $u^k$, we further get
\begin{equation*}
	f_k = -\frac{1}{s_1^{k+1}} \frac{1-s_1}{1-zs_2}= - z^{k+1}s_2^{k+1} \frac{1-1/(zs_2)}{1-zs_2}
	=- z^{k}s_2^{k} \frac{zs_2-1}{1-zs_2}=z^{k}s_2^{k} .
\end{equation*}
Since $G(u)=\dfrac{F(u)-1-zuF(u)}{zu}$, we also find
\begin{equation*}
g_k=\frac1zf_{k+1}-f_k=z^k(s_2^{k+1}-s_2^{k}).
\end{equation*}
We can also compute $\textsc{total}(z)=F(1,z)+G(1,z)$ which counts path that end anywhere, and the result is
\begin{equation*}
\textsc{total}(z)= {\frac {1-z-{z}^{2}-\sqrt { -{z}^{2}-2{z}^{3}-2z+{z}^{4}+1}}{2{z}^{3}}}=\frac1{z^2}g_0.
\end{equation*}
 In retrospective, this is not surprising, since if we consider paths that end at state 0 in the bottom layer, and we go back the last 2 steps, we could have been indeed in any state.
 
 It is worthwhile to notice that 
 \begin{equation*}
f_0+g_0=1+z^2+z^3+2z^4+4z^5+8z^6+17z^7+37z^8+82z^9+185z^{10}+423z^{11}+\cdots 
 \end{equation*}
and the coefficients $1,1,2,4,8,17,\dots$ are sequence A004148 in \cite{OEIS}.

\begin{theorem} The generating functions describing partial Dyck paths with air pockets, landing in state $k$ of the upper/lower layer, are given by
	\begin{equation*}
		f_k=z^{k}s_2^{k},\quad g_k =z^k(s_2^{k+1}-s_2^{k}).
	\end{equation*}
In particular, $f_k+g_k=z^ks_2^{k+1}$ is the generating function of partial paths ending at level $k$.
	\end{theorem}

\section{Right to left model}

Reading Dyck paths with air pockets from right to left means to have arbitrary long up-steps, but only one at the time. While the enumeration for those paths that end at the $x$-axis is the same as before, this is not the case for \emph{partial} paths. 

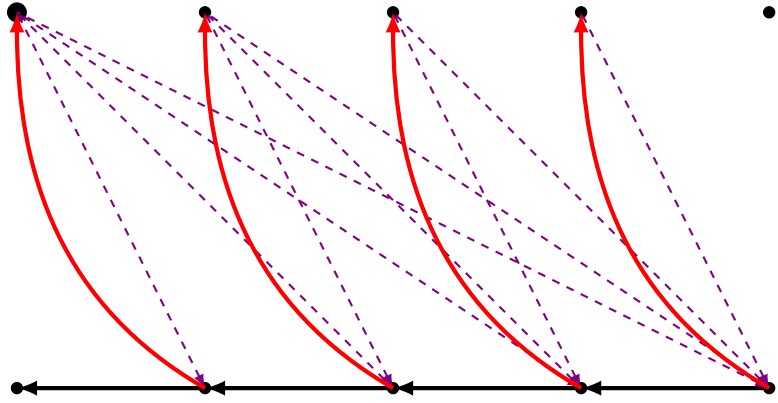
\begin{figure}[ht]\label{myfig2}

	\begin{center}
		\begin{tikzpicture}[scale=2.5]
			\draw (0,2) circle (0.05cm);
			\fill (0,2) circle (0.05cm);
			
			\foreach \x in {0,1,2,3,4}
			{
				\draw (\x,0) circle (0.03cm);
				\fill (\x,0) circle (0.03cm);
			}

			\foreach \x in {0,1,2,3,4}
			{
				\draw (\x,2) circle (0.03cm);
				\fill (\x,2) circle (0.03cm);
			}
			
			\foreach \x in {0,1,2,3}
			{
				\draw[ ultra thick, latex-] (\x,0) -- (\x+1,0);
				
			}

			\draw[   thick, violet,dashed, latex-] (4,0) -- (3,2);
			\draw[   thick, violet,dashed, latex-] (4,0) -- (2,2);
			\draw[   thick, violet,dashed, latex-] (4,0) -- (1,2);
			\draw[   thick, violet,dashed, latex-] (4,0) -- (0,2);
			
			\draw[   thick, violet,dashed, latex-] (3,0) -- (2,2);
			\draw[   thick, violet,dashed, latex-] (3,0) -- (1,2);
			\draw[   thick, violet,dashed, latex-] (3,0) -- (0,2);
			\draw[   thick, violet,dashed, latex-] (2,0) -- (1,2);
			\draw[   thick, violet,dashed, latex-] (2,0) -- (0,2);
			\draw[   thick, violet,dashed,latex-] (1,0) -- (0,2);

			\draw[ ultra   thick, red , latex-] (0,2) to [out= -90,in= 150] (1,0);
			\draw[ ultra   thick, red , latex-] (1,2) to [out= -90,in= 150] (2,0);
			\draw[ ultra   thick, red , latex-] (2,2) to [out= -90,in= 150] (3,0);
			\draw[ ultra   thick, red , latex-] (3,2) to [out=- 90,in= 150] (4,0);

		\end{tikzpicture}
	\end{center}
	\caption{Graphical description of Dyck paths with air pockets. Top layer describes the situation after a down-step, bottom layer after an up-step. }
	
\end{figure}

Figure \ref{myfig2} explains the concept. The generating functions $a_k$ refer to the top layer and $b_k$ to the bottom layer.

The recursions are\footnote{Iverson's notation is used here.} 
\begin{align*}
	a_k&=[k=0]+zb_{k+1},\\
	b_k&=zb_{k+1}+z\sum_{0\le j<k}a_j.
\end{align*}
With bivariate generating functions analogously to before, we find by summing
\begin{equation*}
	A(u)=1+\frac{z}{u}(B(u)-b_0)
\end{equation*}
and
\begin{equation*}
	B(u)=\frac{z}{u}(B(u)-b_0)+z\sum_{0\le j<k}a_ju^k=\frac{z}{u}(B(u)-b_0)+\frac{zu}{1-u}A(u).
\end{equation*}
One variable can be eliminated:
\begin{equation*}
	B(u)= \frac{z}{u}(B(u)-b_0)+\frac{zu}{1-u} +\frac{z^2}{1-u}(B(u)-b_0) .
\end{equation*}
Solving
\begin{equation*}
	B(u)= {\frac {z \left( B(0)-B(0)u-{u}^{2}+zB(0)u \right) }{z
			-zu+{z}^{2}u-u+{u}^{2}}}
		\end{equation*}
 The denominator factors as $(u-s_1^{-1})(u-s_2^{-1})$.
The bad factor is this time $(u-s_1^{-1})$. Dividing it out,
\begin{equation*}
	B(u)=\frac {z \left( -us_1-B(0)s_1+B(0)s_1
		z-1 \right) }{us_1-zs_1+{z}^{2}s_1-s_1+1}
		\end{equation*}
and further
\begin{equation*}
	B(0)= b_0=\frac{z}{s_1-1}=s_2-1.
\end{equation*}
Thus, after some simplifications,
 \begin{equation*}
	B(u)=-z+\frac{zs_1}{(s_1-1)(1-\frac{u}{zs_1})},
\end{equation*}
or
\begin{equation*}
	B(u)=-z+\frac{1}{s_2(s_1-1)(1-s_2u)}=
	-z+\frac{s_2-1}{zs_2 (1-s_2u)}
\end{equation*}
and then
\begin{equation*}
b_k=\frac{s_2-1}{z  }s_2^{k-1},\quad k\ge1.
\end{equation*}
The functions $a_k$ could be computed from here as well, but for the partial paths only the functions $b_k$ are of relevance, if we don't consider the empty path.

\begin{theorem}
	The generating functions of partial Dyck paths with air pockets in the right to left model are
	\begin{align*}
1+b_0=s_2
	\end{align*}
and
\begin{equation*}
	b_k=\frac{s_2-1}{z  }s_2^{k-1},\quad k\ge1.
\end{equation*}
To consider the total does not make sense, since in just 1 or 2 steps, every state can be reached, so a sum over $b_k$ would not converge.
\end{theorem}

\section{Skew Dyck paths with air pockets}

\begin{figure}[h]

	\begin{center}
		\begin{tikzpicture}[scale=1.5]
			\draw (0,0) circle (0.1cm);
			\fill (0,0) circle (0.1cm);
			
			\foreach \x in {0,1,2,3,4,5,6,7,8}
			{
				\draw (\x,0) circle (0.05cm);
				\fill (\x,0) circle (0.05cm);
			}
			
			\foreach \x in {0,1,2,3,4,5,6,7,8}
			{
				\draw (\x,-1) circle (0.05cm);
				\fill (\x,-1) circle (0.05cm);
			}
			
			\foreach \x in {0,1,2,3,4,5,6,7,8}
			{
				\draw (\x,-2) circle (0.05cm);
				\fill (\x,-2) circle (0.05cm);
			}
			
			\foreach \x in {0,1,2,3,4,5,6,7}
			{
				\draw[ thick,-latex] (\x,0) -- (\x+1,0);
				
			}

			\foreach \x in {1,2,3,4,5,6,7}
			{
				\draw[thick,  latex-] (\x+1,0) to[out=200,in=70]  (\x,-1);

			}
			\draw[ thick,  latex-] (1,0) to[out=200,in=70]  (0,-1);

			\foreach \x in {0,1,2,3,4,5,6,7}
			{
				
				\draw[thick,  latex-] (\x,-1) to[out=30,in=250]  (\x+1,0);	
				
			}

			\foreach \x in {0,1,2,3,4,5,6,7}
			{
				\draw[ thick,-latex] (\x+1,-1) -- (\x,-1);
				
			}
			\foreach \x in {0,1,2,3,4,5,6,7}
			{
				\draw[ thick,-latex,red] (\x+1,-1) -- (\x,-2);
				
			}
			
			\foreach \x in {0,1,2,3,4,5,6,7}
			{
				\draw[ thick,-latex,red] (\x+1,-2) -- (\x,-2);
				
			}
			
			\foreach \x in {0,1,2,3,4,5,6,7}
			{
				\draw[ thick,-latex] (\x+1,-2) -- (\x,-1);
				
			}

		\end{tikzpicture}
	\end{center}
	\caption{Three layers of states according to the type of steps leading to them (up, down-black, down-red). }
	\label{threelayers}
\end{figure}
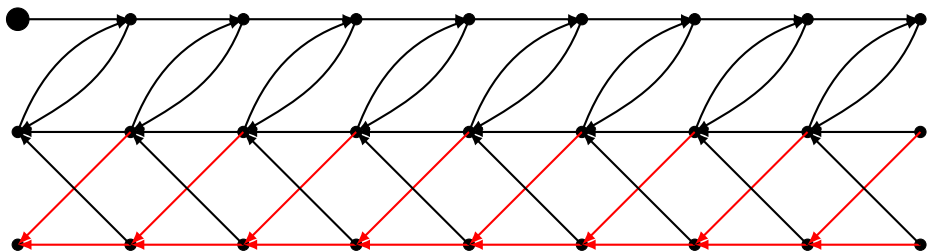

The walks according to Figure~\ref{threelayers} are related to skew Dyck paths  \cite{skew-paper}; the red down-steps are modeled to stand for south-west steps, 
and the way they are arranged, there are no overlaps of such a path. See \cite{skew-paper} and the references cited there.

\begin{figure}[h]

	\begin{center}
		\begin{tikzpicture}[scale=2.5]
			\draw (0,0) circle (0.07cm);
			\fill (0,0) circle (0.07cm);
			
			\foreach \x in {0,1,2,3,4}
			{
				\draw (\x,0) circle (0.05cm);
				\fill (\x,0) circle (0.05cm);
			}
			
			\foreach \x in {0,1,2,3,4 }
			{
				\draw (\x,-1) circle (0.05cm);
				\fill (\x,-1) circle (0.05cm);
			}
			
			\foreach \x in {0,1,2,3,4 }
			{
				\draw (\x,-2) circle (0.05cm);
				\fill (\x,-2) circle (0.05cm);
			}
			
			\foreach \x in {0,1,2,3  }
			{
				\draw[ultra thick,-latex] (\x,0) -- (\x+1,0);
				
			}

			\foreach \x in {0,1,2,3 }
			{

				\draw[ ultra thick,  -latex] (\x,-1) to[out=90,in=190]  (\x+1,0);	
				
			}
			
			\draw[   thick, gray,dashed, -latex] (4,0) -- (3,-1);
			\draw[   thick, gray,dashed, -latex] (4,0) -- (2,-1);
			\draw[   thick, gray,dashed, -latex] (4,0) -- (1,-1);
			\draw[   thick, gray,dashed, -latex] (4,0) -- (0,-1);
			\draw[   thick, gray,dashed, -latex] (3,0) -- (2,-1);
			\draw[   thick, gray,dashed, -latex] (3,0) -- (1,-1);
			\draw[   thick, gray,dashed, -latex] (3,0) -- (0,-1);
			\draw[   thick, gray,dashed, -latex] (2,0) -- (1,-1);
			\draw[   thick, gray,dashed, -latex] (2,0) -- (0,-1);
			\draw[   thick,  gray,dashed, -latex] (1,0) -- (0,-1);

			\draw[ ultra  thick, red, -latex] (4,-1) -- (3,-2);
			\draw[  ultra thick,  red, -latex] (3,-1) -- (2,-2);
			\draw[  ultra thick,  red, -latex] (2,-1) -- (1,-2);

			\draw[  ultra thick,  red, -latex] (1,-1) -- (0,-2);
				\draw[  ultra thick,  red, -latex] (2,-2) -- (1,-2);
								\draw[  ultra thick,  red, -latex] (3,-2) -- (2,-2);
												\draw[  ultra thick,  red, -latex] (4,-2) -- (3,-2);
	\draw[  ultra thick,  red, -latex] (1,-2) -- (0,-2);
			
			\draw[   thick, gray,dashed, -latex] (4,-2) -- (3,-1);
			\draw[   thick,  gray,dashed, -latex] (4,-2) -- (2,-1);
			\draw[   thick,  gray,dashed, -latex] (4,-2) -- (1,-1);
			\draw[   thick,  gray,dashed, -latex] (4,-2) -- (0,-1);
			\draw[   thick,  gray,dashed, -latex] (3,-2) -- (2,-1);
			\draw[   thick,  gray,dashed, -latex] (3,-2) -- (1,-1);
			\draw[   thick,  gray,dashed, -latex] (3,-2) -- (0,-1);
			\draw[   thick,  gray,dashed, -latex] (2,-2) -- (1,-1);
			\draw[   thick,  gray,dashed, -latex] (2,-2) -- (0,-1);
			\draw[   thick,  gray,dashed, -latex] (1,-2) -- (0,-1);

			\foreach \x in {0,1,2,3,4 }
			{
				
			}

		\end{tikzpicture}
	\end{center}
	\caption{Three layers of states according to the type of steps leading to them (up, down-black, down-red). Black down-steps are condensed into giant grey down-steps.}
	\label{threelayers2}
\end{figure}
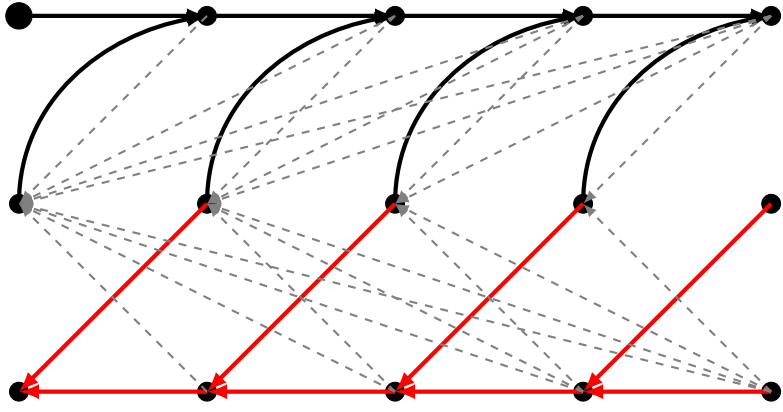

Now we combine this model with air pockets. Each maximal sequence of black down-steps is condensed into one giant down-step, depicted
in dashed grey in Figure \ref{threelayers2}

Introducing generating functions, according to the three layers, we find the following recursions by inspection;
 \begin{align*}
	a_0&=1,\quad a_{k+1}=za_k+zb_k,\ k\ge 0,\\
	b_k&=z\sum_{j>k}a_j+z\sum_{j>k}c_j,\\
	c_k&= zb_{k+1}+zc_{k+1}.
\end{align*}
Translating these into bivariate generating functions, we further have
\begin{gather*}
A(u)=1+zuA(u)+zuB(u),\\
B(u)=\frac z{1-u}[A(1)-A(u)]+\frac z{1-u}[C(1)-C(u)],\\
	C(u)=zuB(u)+zuC(u).
\end{gather*} 
Solving,
\begin{align*}
A(u)&={\frac {{z}^{3}{u}^{2}A(1)+{z}^{3}{u}^{2}C(1)-z{u}^{2}-{z}^{2}
		uC(1)-{z}^{2}u-{z}^{2}uA(1)+zu+u-1}{ \left( -1+zu \right) 
		\left( z{u}^{2}+2 {z}^{2}u-zu-u+1 \right) }},\\
	B(u)&=-{\frac { \left( -A(1)-C(1)+zuA(1)+zuC(1)+1 \right) z}
		{z{u}^{2}+2 {z}^{2}u-zu-u+1}},\\
		C(u)&={\frac {{z}^{2}u \left( -A(1)-C(1)+zuA(1)+zuC(1)+1
				\right) }{ \left( -1+zu \right)  \left( z{u}^{2}+2 {z}^{2}u-zu-u+1
				\right) }}.
		\end{align*}
	We factor $z{u}^{2}+2 {z}^{2}u-zu-u+1=(u-s_1)(u-s_2)$
	with	
	\begin{equation*}
s_2={\frac {-2{z}^{2}+z+1-\sqrt {4 {z}^{4}-4 {z}^{3}-3 {z}^{2}-2			 z+1}}{2z}},\quad s_1=\frac{1}{zs_2}.
	\end{equation*}
Since $A(u)-C(u)=\frac{1}{1-zu}$, we have $A(1)-C(1)=\frac{1}{1-z}$, and we only need to compute one of them.
Dividing the (bad) factor $(u-s_2)$ out, plugging in $u=1$ and solving leads to 
\begin{equation*}
A(1)={\frac {-s_2z+2-z}{ 2\left(1- s_2z \right)  \left(1- z \right) }}=\frac1{2(1-z)}+\frac1{2(1-zs_2)}
\end{equation*}
and
\begin{equation*}
C(1)=-\frac1{2(1-z)}+\frac1{2(1-zs_2)}.
\end{equation*}
Using these values, we find
\begin{equation*}
A(u)+B(u)+C(u)=\frac{s_2(1-z^2-zs_2)}{(1-zs_2)(1-uzs_2)}
\end{equation*}
and furthermore
\begin{equation*}
	[u^k](A(u)+B(u)+C(u))=\frac{z^ks_2^{k+1}(1-z^2-zs_2)}{(1-zs_2)}
\end{equation*}
These functions describe all skew Dyck paths with air pockets, ending at level $k$.
For $k=0$, this yields
\begin{equation*}
1+z^2+z^3+3z^4+7z^5+17z^6+45z^7+119z^8+323z^9+893z^{10}+2497z^{11}+\cdots.
\end{equation*}


\begin{thebibliography}{99}
	
	
	
	\bibitem{valentine}
	Jean-Luc Baril, Sergey Kirgizov, Rémi Maréchal, Vincent Vajnovszki,
		Enumeration of Dyck paths with air pockets,
	arXiv:2202.06893.
	
	
	\bibitem{prodinger-kernel}
	H. Prodinger,
	\newblock { The kernel method: A collection of examples}, 
	\newblock {\em S\'em. Lothar. Combin.},  {\bf B50f} (2004), 19 pages.
	
	
	\bibitem{skew-paper} H. Prodinger, Partial skew Dyck paths---a kernel method approach, preprint, 2021. 
  
	
	
	
	\bibitem{OEIS}
	N. J. A.~Sloane et al., The On-line Encyclopedia of Integer Sequences, 2022. 
	
		\end{thebibliography}
\end{document}